\documentclass{amsart}
	\usepackage{array}
	\usepackage{amsmath}
	\usepackage{amsfonts,amssymb}
\usepackage{amsthm}

\newcolumntype{C}{>{$}c<{$}}

\newtheorem{theorem}{Theorem}[section]
\newtheorem{definition}{Definition}[section]
\newtheorem{proposition}{Proposition}[section]
\newtheorem{remark}{Remark}[section]
\newtheorem{lemma}{Lemma}[section]

\newtheorem{corollary}{Corollary}[section]

\def\gc{\mathfrak g_{{}_{\mathbf C}}}
\def\kc{\mathfrak k_{{}_{\mathbf C}}}
\def\pc{\mathfrak p_{{}_{\mathbf C}}}

\def\Kc{K_{{}_{\mathbf C}}}

\def\Gc{G_{{}_{\mathbf C}}}

\def\g{\mathfrak g}
\def\k{\mathfrak k}
\def\s{\mathfrak s}
\def\p{\mathfrak p}
\def\a{\mathfrak a}

\def\S{{\mathbf S}_{\mathbf R}}
\def\sr{{\mathbf s}_{\mathbf R}}

\begin{document}
\title{Complexity of nilpotent orbits and the Kostant-Sekiguchi correspondence}
\author{Donald R. King}
\begin{abstract}
Let $G$ be a connected linear semisimple Lie group with Lie algebra $\mathfrak g$, and let ${\Kc}~\rightarrow~{Aut(\pc)}$ be the complexified isotropy representation at the identity coset of the corresponding symmetric space. Suppose that $\Omega$ is a nilpotent $G$-orbit in $\mathfrak g$ and $\mathcal O$ is the nilpotent $\Kc$-orbit in $\pc$ associated to $\Omega$ by the Kostant-Sekiguchi correspondence. We show that the complexity of $\mathcal O$ as a $\Kc$ variety measures the failure of the Poisson algebra of smooth $K$-invariant functions on $\Omega$ to be commutative. 
\end{abstract}
\maketitle
\markboth{D. R. King}{Complexity of Nilpotent Orbits}

\section{Introduction}

	The Kostant-Sekiguchi correspondence is a vital tool in the study of infinite dimensional representations of semisimple Lie groups. Let us recall some facts about this correspondence in case $G$ is a connected real, linear semisimple Lie group with maximal compact subgroup $K$. We obtain the corresponding Cartan decomposition $\g=\k \oplus\p$ where $\g$ (resp. $\k$) is the Lie algebra of $G$ (resp., $K$). The vector spaces $\g$, $\k$ and $\p$ are then complexified to give a vector space decomposition of $\gc$,  the Lie algebra of $\Gc$ (the complexification of $G$), as $\gc=\kc\oplus \pc$. The Kostant-Sekiguchi correspondence is a bijection between the nilpotent $G$ orbits in $\g$ and the nilpotent $\Kc$ orbits in $\pc$. (For the precise definition of the correspondence we refer the reader to \cite{king}.) \par
	If $\Omega$ is a nilpotent $G$-orbit in $\mathfrak g$ and $\mathcal O$ is the nilpotent $\Kc$-orbit in $\pc$ associated to $\Omega$ by the Kostant-Sekiguchi correspondence, then $(\Omega,\ \mathcal O)$ is said to be a Kostant-Sekiguchi pair. Among the nice elementary properties of such a pair are: (1) $\Omega$ and $\mathcal O$ lie in the same $\Gc$ orbit which we denote by  ${\mathcal O}_{\mathbf C}$ and (2) $\mathcal O$ is a Lagrangian submanifold of ${\mathcal O}_{\mathbf C}$ (relative to the Kostant-Souriau symplectic form on ${\mathcal O}_{\mathbf C}$.) Moreover, Vergne \cite{ver} has established a much deeper relationship between $\Omega$ and $\mathcal O$, namely that there is a $K$-equivariant diffeomorphism which maps $\Omega$ onto $\mathcal O$.  \par
	Recently, the author proved that if $(\Omega,\ \mathcal O)$ form a Kostant-Sekiguchi pair then $\Omega$ is multiplicity free as a Hamiltonian $K$-space if and only if $\mathcal O$ is a spherical $\Kc$ variety \cite{king}. (The definition of multiplicity free is given below in Remark \ref{multfree}. Spherical $\Kc$-varieties are defined below in Remark \ref{spherical}.) The goal of this paper is to prove a generalization of that result to all Kostant-Sekiguchi pairs. That generalization is contained in Theorem \ref{maintheorem}. In essence our theorem shows that the complexity of the $\Kc$ action on $\mathcal O$ measures the failure of the Poisson algebra of smooth $K$-invariant functions on $\Omega$ to be commutative. \par
	The author wishes to thank A. T. Huckleberry for his proof of the density result in Proposition \ref{keydensityresult} and Maxim Braverman for some useful discussions.

\section{Notation and key definitions}
	We now introduce some further concepts and notations. Unless otherwise indicated, in this section $K$ will denote an arbitrary compact group. We assume that $K$ is contained in its complexification $\Kc$. Our basic reference for symplectic manifolds with Hamiltonian $K$-actions is \cite{gust2}. \par

	If $A\subset\Gc$ is a Lie subgroup and $S\subset \gc$ then $A^S$  denotes the subgroup of $A$ that fixes each element of $S$ under the adjoint action of $\Gc$ on $\gc$. If $\a\subset\gc$ is a Lie subalgebra, then $\a^S$ is defined similarly.\par

For the remainder of this section $\mathbf X$ will denote a connected symplectic manifold which is a Hamiltonian $K$-space. Let the $K$-invariant symplectic form be $w_{\mathbf X}$ and the moment map be $\Phi:\mathbf X\rightarrow \k^*$. 
\begin{definition}\label{collective}
Let $\mathfrak A = C^{\infty}(\mathbf X)^K$ be the algebra (with respect to Poisson bracket) of $K$ invariant smooth functions on $\mathbf X$. $\mathfrak Z$ denotes the center of $\mathfrak A$. 
\end{definition} 
\begin{remark}\label{multfree} $\mathbf X$ is said to be multiplicity free if  $\mathfrak A=\mathfrak Z$  i.e., $\mathfrak A$ is a commutative Poisson algebra. 
\end{remark}

	In addition, we need to recall some facts about the symplectic structure of $\Omega$ and the action of $K$ on $\Omega$. ($K$ is a maximal compact subgroup of $G$.) Suppose that $\Omega=G\cdot E$. For each $E'\in\Omega$, we identify $T_{E'}(\Omega)$, the tangent space of $\Omega$ at $E'$, with the quotient $\g\slash {\g}^{E'}$. The Kostant-Souriau form $w_{\Omega}$ on $\Omega$ is defined by setting 
\begin{equation*}
w_{\Omega}|_{E'}(\bar{Y},\ \bar{Z})= \kappa(E',\ [Y,\ Z])
\end{equation*}
, for all $Y,\ Z\in\g$ where $\kappa$ denotes the Killing form of $\g$. ($\bar{Y}$ and $\bar{Z}$ are the cosets of $Y$ and $Z$ in $\g\slash {\g}^{E'}$.) $w_{\Omega}$ is a symplectic form on $\Omega$.\par 
	The action of $K$ on $\Omega$ is the left action. If $\xi\in\mathfrak{k}$, then $\xi$ determines a global vector field $\xi_{\Omega}$ on ${\Omega}$ according to the following definition:
$$ \xi_{\Omega}f(E') = {\frac{df(exp(-t\xi)\cdot E')}{dt}\Bigg\vert}_{t=0}$$
where $f$ is any smooth function on ${\Omega}$ and $E'\in {\Omega}$. If $Y$ is any smooth vector field on $\Omega$, $Y(E')$ is the vector in $T_{E'}(\Omega)$ obtained by evaluating $Y$ at $E'$. \par
	Since $w_{\Omega}$ is invariant under $G$, it is invariant under $K$. The $K$ action on $\Omega$ is Hamiltonian in the following sense. For each $\xi\in \mathfrak{k}$, there is a function $\phi^{\xi}\in C^{\infty}({\Omega})$ such that $\iota(\xi_{\Omega})w_{\Omega} = -d\phi^{\xi}$, where $\iota(\xi_{\Omega})$ denotes interior multiplication by $\xi_{\Omega}$. We may take $\phi^{\xi}=\kappa(\xi,\ \cdot)$. We obtain the moment mapping $\Phi_{\Omega}: {\Omega}\rightarrow {\mathfrak k}^*$ by setting $\Phi_{\Omega}(E')(\xi) = \phi^{\xi}(E')=\kappa(E',\ \xi)=\kappa(E_k',\ \xi)$ for all $E'\in {\Omega}$ and $\xi\in\mathfrak{k}$. $E'_k$ denotes the component of $E'$ in $\k$. 
\par 

	Thus, $\Omega$ is a Hamiltonian $K$-space, with moment mapping $\Phi_{\Omega}:\Omega\rightarrow \k$ defined by sending an element $E'\in \Omega$ to  its component in $\k$. (We have identified $\k$ with its real dual space $\k^*$ using the restriction of $\kappa$ to $\k$.) Each function $f$ in $C^{\infty}({\Omega})$ gives rise to a smooth vector field $X_f$ satisfying $df(Y)=w_{\Omega}(Y,\ X_f)$ for all smooth vector fields $Y$ on $\Omega$. $X_f$ is said to be the Hamiltonian vector field associated to $f$. (If $\xi\in\k$, then $X_{\phi^{\xi}}=\xi_{\Omega}$.)  $C^{\infty}({\Omega})$ is a Poisson algebra under the Poisson bracket $\{\cdot,\ \cdot\}$ defined as follows: for $f,\ g\in C^{\infty}({\Omega})$,  $\{f,\ g\}=w_{\Omega}(X_f,\ X_g)$. One checks that the linear mapping $\xi\mapsto {\phi}^{\xi}$ is a Lie algebra homomorphism $\k\rightarrow C^{\infty}({\Omega})$ when $\{\cdot,\ \cdot\}$ is taken as the Lie bracket on $C^{\infty}({\Omega})$.\par

\begin{lemma}\label{tanspacedecomp} Let $E'\in\mathbf X$. Set $W=W_{E'}= T_{E'}(K\cdot E')$ and $W^{\bot}$ equal to the orthogonal complement (with respect to $\zeta=w_{\mathbf X}|_{E'}$) of $W$ inside $T_{E'}(\mathbf X)$, then we have the following orthogonal decomposition of $T_{E'}(\mathbf X)$ into (real) symplectic vector spaces:
\begin{equation*}
T_{E'}(\mathbf X)= \frac{W}{W\cap W^{\bot}}\bigoplus\left( (W\cap W^{\bot})\oplus (W\cap W^{\bot})^*\right)\bigoplus \frac{W^{\bot}}{W\cap W^{\bot}}.
\end{equation*}
The restriction of $\zeta$ to $(W\cap W^{\bot})\oplus (W\cap W^{\bot})^*$ is given by $\zeta((Y_1,\ \lambda_1),\ (Y_2,\ \lambda_2))=\lambda_1(Y_2) -\lambda_2(Y_1)$ for all $Y_i\in (W\cap W^{\bot})$, $\lambda_i\in(W\cap W^{\bot})^*$.  Moreover, (1) ${{\k}^{{\Phi}(E')}\slash {\k}^{E'}}$ and ${W\cap W^{\bot}}$ are isomorphic (as $\k^{E'}$ modules) and (2) ${\k}\slash{{\k}^{{\Phi}(E')}}$ and $\displaystyle{\frac{W}{W\cap W^{\bot}}}$ are isomorphic (as $\k^{E'}$ modules).

\end{lemma}
\begin{proof} See Corollary 9.10 and Lemma 9.11 of \cite{karshon}. These results are steps in the proof of the normal slice theorem of Guillemin, Sternberg and Marle. \end{proof}
\begin{definition}\label{codimomega} If $E''\in {\mathbf X}$ (or $\k$), set $d(E'')=\dim K\cdot E''$. Let $d=d_{{\mathbf X}}$ be the maximum dimension of a $K$-orbit in ${\mathbf X}$ and let $m=m_{{\mathbf X}}$ denote the minimal codimension of a $K$-orbit in ${\mathbf X}$. $d_{\Phi}$ will denote the maximum dimension of a $K$ orbit in $\Phi({\mathbf X})$. We define three important subsets of ${\mathbf X}$.
\begin{equation}\label{Omegasubsets}
\begin{split}
&{\mathbf X}_d=\{E'\in {\mathbf X}|\text{$\dim {K\cdot E'}=d$}\}\cr
&{\mathbf X}_0=\{E'\in {\mathbf X}|\text{$\exists$ an open set $U\subset\mathbf X$ such that $E'\in U$ and $d(\cdot)$ is constant on $U$}\}\cr
&{\mathbf X}_{\Phi}=\{E'\in {\mathbf X}|\text{$\dim {K\cdot \Phi(E')}=d_{\Phi}$}\}=\{E'\in {\mathbf X}|\text{$\dim {\k}^{ \Phi(E')}$ is minimum}\}
\end{split}
\end{equation}
\end{definition}

	\begin{remark}\label{Omegaprops}
	 From Proposition 27.1 in \cite{gust2} we know that ${\mathbf X}_0$ is open and dense in ${\mathbf X}$. In addition, ${\mathbf X}_0\subset {\mathbf X}_d$. Otherwise there is a point in ${\mathbf X}_0$ whose $K$ orbit has dimension $d'<d$. But, then there must be an open subset of points whose $K$ orbits have dimension $d'$ which is impossible. \par
\end{remark}
\begin{proposition}\label{keydensityresult}
${\mathbf X}_{\Phi}$ is open and dense in ${\mathbf X}$. 
\end{proposition}
\begin{proof} The fact that ${\mathbf X}_{\Phi}$ is open is Lemma 1 in section 3 of \cite{huckwurz}. The following proof (unpublished) that ${\mathbf X}_{\Phi}$ is dense in $\mathbf X$ is due to A. T. Huckleberry.\par
	It suffices to show that ${\mathbf X}_{\Phi}^c$, the complement of ${\mathbf X}_{\Phi}$, has codimension at least 2 in $\mathbf X$. We use induction on the dimension of $K$. If $\dim K=0$, then $K$ is finite. Then ${\mathbf X}_{\Phi} = \mathbf X$ since all $K$ orbits in $\Phi(\mathbf X)$ are zero dimensional. Suppose $\dim K>0$. If $x\in {\mathbf X}_{\Phi}^c$, we need to show that there is an open neighborhood of $x$ whose intersection with ${\mathbf X}_{\Phi}^c$ has codimension at least 2. Let $K_0$ denote the identity component of $K$. First consider ${\mathbf X}^{K_0}$, the set of fixed points of $K_0$ on $\mathbf X$, and its intersection with ${\mathbf X}_{\Phi}^c$.  If $K_0\cdot x=x$, we use the slice theorem to construct an open neighborhood $U=K\times_{K^x} \Sigma$ of $x$. Define $\Sigma_{\Phi}$ relative to the action of $({K^x})_0$ in the obvious way. The argument for Proposition 27.3 in \cite{gust2} shows that either (a) $K_0$ (which equals $({K^x})_0$) acts trivially on $\Sigma$, or (b) $U^{K_0}$ has codimension $\ge$ 2 in $U$. In case (a) $K$ acts as a finite group on $U$ so that $U=U_{\Phi}$. In case (b) $U^{K_0}\cap (U_{\Phi})^c$ has codimension $\ge$ 2 in $U$. This argument takes care of the points in ${\mathbf X}_{\Phi}^c\cap {\mathbf X}^{K_0}$. If $x\notin {\mathbf X}^{K_0}$, then the slice neighborhood $U = K\times_{K^x} \Sigma$ has the property that $\dim K^x<\dim K$. By induction $(\Sigma_{\Phi})^c$ has codimension at least 2 in $\Sigma$. Thus $(U_{\Phi})^c=K\times_{K^x} ({\Sigma}_{\Phi})^c$ has codimension at least 2 in $U$.
\end{proof}
\par
Since ${\mathbf X}_d$ and ${\mathbf X}_{\Phi}$ are each open and dense in ${\mathbf X}$, we have the following result.
\begin{corollary} ${\mathbf X}_{\Phi}\cap {\mathbf X}_d$ is open and dense in ${\mathbf X}$.
\end{corollary}
We know recall the notion of a coisotropic submanifold of ${\mathbf X}$.
\begin{definition}
The orbit $K\cdot E'$ in ${\mathbf X}$ is said to be coisotropic if $W_{E'}^{\bot}\subset W_{E'}$.
\end{definition}
	By results of Guillemin and Sternberg \cite{gust1}, ${\mathbf X}$ is multiplicity free as a Hamiltonian $K$-space if and only if there is an open dense subset $U$ of ${\mathbf X}$ such that for $E'\in U$ the orbit $K\cdot E'$ is coisotropic in ${\mathbf X}$. Therefore, it is reasonable to use the size of the quotient $\displaystyle {\frac{{T_{E'}(K\cdot E')}^{\bot}}{{T_{E'}(K\cdot E')}\cap {T_{E'}(K\cdot E')}^{\bot}}}$ for generic $K$ orbits to measure the failure of ${\mathbf X}$ to be multiplicity free. 
\begin{proposition}\label{epsilonofomega}
 There is an open dense subset $U$ of ${\mathbf X}$ such that for all $E'\in U$, the non-negative integer 
\begin{equation}\label{epsilonequality}
\dim{\frac{{T_{E'}(K\cdot E')}^{\bot}}{{T_{E'}(K\cdot E')}\cap {T_{E'}(K\cdot E')}^{\bot}}},
\end{equation}
has a constant value. We denote this value by $2\epsilon({\mathbf X})$.
\end{proposition}

\begin{proof}
Choose $U={\mathbf X}_d\cap {\mathbf X}_{\Phi}$. Then for all $E'\in U$, the decomposition of $T_{E'}({\mathbf X})$ in Lemma \ref{tanspacedecomp} implies that 
\begin{equation}\label{formulaforepsilon}
\dim{\frac{{T_{E'}(K\cdot E')}^{\bot}}{{T_{E'}(K\cdot E')}\cap {T_{E'}(K\cdot E')}^{\bot}}}=\dim {\mathbf X} - 2d +d_{\Phi}
\end{equation}

Note that if there is another dense open set $U'$ on which the dimension function in (\ref{epsilonequality}) is constant, then $U'$ must have non-empty intersection with ${\mathbf X}_d\cap {\mathbf X}_{\Phi}$. Thus,  $\epsilon({\mathbf X})$ is well defined.
\end{proof}
As a corollary of the preceding proof we have
\begin{corollary}\label{symsliceeq} Suppose that $E'\in {\mathbf X}_d\cap {\mathbf X}_{\Phi}$, then
\begin{equation*}
\dim{\frac{{T_{E'}(K\cdot E')}^{\bot}}{{T_{E'}(K\cdot E')}\cap {T_{E'}(K\cdot E')}^{\bot}}}=2\epsilon({\mathbf X}).
\end{equation*}
\end{corollary}
	We recall from \cite{huckwurz} or \cite{king} the definition of the rank of the action of $K$ on $X$.
\begin{definition}
The rank of the $K$-action on $X$, denoted by $r_K(X)$, is equal to rank $K$ - rank $K^{E'}$, where $E'\in X$ and the orbit $K\cdot E'$ has maximum dimension among the $K$ orbits in $X$.
\end{definition}
\medskip
	Finally, we recall from \cite{pa1}, the notions of rank and complexity of algebraic $\Kc$ actions.
\begin{definition}\label{complexity} Suppose that $Y$ is a variety with $\Kc$ action and $B_k$ is a Borel subgroup of $\Kc$. The complex codimension of a generic $B_k$ orbit is called the complexity of $Y$, denoted $c_{\Kc}(Y)$ or $c(Y)$ (when the reductive group $\Kc$ is understood).  If $U_k$ is the nilpotent radical of $B_k$ and $B_k\cdot z$ is a generic $B_k$ orbit, then the codimension of $U_k\cdot z$ in $B_k\cdot z$ is called the rank of $Y$. It is denoted $r_{\Kc}(Y)$. 
\end{definition}
\begin{remark}\label{spherical} $c(Y)$ is also the transcendence degree (over $\mathbf C$) of the $B_k$ invariant functions in the field of rational functions (with complex coefficients) on $Y$. $Y$ is spherical for $\Kc$ if and only if $c(Y)=0$. The rank of $Y$ is also the transcendence degree of the $U_k$ invariants in the field of rational functions on $Y$.
\end{remark}
\section{Main Theorem}
	Our main result is: 
\begin{theorem}\label{maintheorem} Let $(\mathcal O,\ \Omega)$ be a Kostant-Sekiguchi pair, then\par  (a) $r_{\Kc}(\mathcal O) = r_K(\Omega)$;\par   (b) $c(\mathcal O)=\epsilon(\Omega)$.
\end{theorem}
	This result was inspired by the main result of \cite{myky1}.\par
We assume from now on that $(\Omega,\ \mathcal O)$ is a Kostant-Sekiguchi pair. The proof of  Theorem \ref{maintheorem} requires two important facts about $\Omega$ and $\mathcal O$. The first is the existence of a $K$ invariant diffeomorphism
\begin{equation}\label{Verdiff}
{\mathcal V}_{\Omega}: \Omega\rightarrow \mathcal O
\end{equation}
 that was mentioned above. The second fact is the existence of a $K$-invariant diffeomorphism 
\begin{equation}\label{conormaldiff}
{\mathcal M}_{\mathcal O}:\mathcal O\rightarrow{K\times_{K^{\displaystyle \s}} V_{\mathcal O}(\s)}
\end{equation}
established in the proof of Proposition 5.2 of \cite{king}. To describe the vector bundle $K\times_{K^{\displaystyle \s}} V_{\mathcal O}(\s)$ in (\ref{conormaldiff}), we recall the notation of \cite{king}.\par
	 There is a Kostant-Sekiguchi $sl(2)$-triple $\{x,\ e, \ f\}$ such that $\mathcal O=\Kc\cdot e$. The Kostant-Sekiguchi property means that (1) $x\in i\k$, $e,\ f\in\pc$, (2) $e=\sigma(f)$, where $\sigma$ is conjugation on $\gc$ relative to the real form $\g$, and (3) the following Lie bracket relations hold: $[x,\ e]=2e$, $[x,\ f]=-2f$, and $[e,\ f]=x$. It follows that the Lie algebra ${\mathbf C}x\oplus{\mathbf C}e\oplus {\mathbf C}f$ is the complexification of a Lie subalgebra $\s$ of $\g$, where $\s$ is isomorphic to $sl(2,\ \mathbf R)$. $V_{\mathcal O}(\s)$ is the quotient $[\kc,\ e]\slash[\k,\ e]$. In \cite{king} it is shown that $K\times_{K^{\displaystyle\s}} V_{\mathcal O}(\s)$ is diffeomorphic to the conormal bundle of $K\cdot e$ inside the cotangent bundle of $\mathcal O$. In addition, there is an isomorphism of $K^{\s}$ modules over $\mathbf R$:
\begin{equation}\label{VsubOofs}
V_{\mathcal O}(\s)\simeq {{\k}^x\slash{{\k}^{\s}}}\oplus Z.
\end{equation}
where $Z$ is the sum (over $\mathbf C$) of the positive eigenspaces of $ad(x)$ on $\kc$ that do not lie in $\kc^e$.\par

We now establish the main result, Theorem \ref{maintheorem}.
\begin{proof}	We first establish part (a).\par

 By composing the $K$-invariant diffeomorphisms in (\ref{Verdiff}) and (\ref{conormaldiff}), we have the assignment: $E'\mapsto {\mathcal M}\circ{\mathcal V}(E')$. So we can identify $E'$ with an equivalence class $[k_0,\ y+z]$ in $K\times_{K^{\displaystyle \s}} V_{\mathcal O}(\s)$ where $y\in{{\k}^x\slash{{\k}^{\s}}}$ and $z\in Z$. (See equation (\ref{VsubOofs}).) It is more convenient to consider the equivalence class ${\mathcal M}\circ{\mathcal V}(k_0^{-1}\cdot E')=[\mathbf{1},\ y+z]$ where $\mathbf{1}$ denotes the identity in $K$. Assume that $\dim K\cdot E'=d$. Since $\dim K\cdot E'$ is maximum, (1) $(\k^{\s})^{y+z}$ has minimum dimension among the subalgebras $(\k^{\s})^{v}$ for $v\in V_{\mathcal O}(\s)$ and (2) $(\k^{\s})^y$ has minimum dimension among the subalgebras $(\k^{\s})^{y'}$ for $y'\in \k^x\slash{{\k}^{\s}}$. (See Lemma 6.2 of \cite{king}.)\par
	 Set ${\mathbf s}_{\mathbf R}=(\k^{\s})^y$ and ${\mathbf s}_{\mathbf C}=(\kc^{\s})^y$. Set $\mathbf S$ equal to the connected subgroup of $\Kc^{\s}$ with Lie algebra ${\mathbf s}_{\mathbf C}$ and $\S$ equal to the connected subgroup of $\mathbf S$ with Lie algebra ${\mathbf s}_{\mathbf R}$. Then $\mathbf S$ (resp., $\S$) is a stabilizer of general position for the action of $\Kc^{\s}$ on $\kc^x\slash{{\kc}^{\s}}$ (resp., $\k^x\slash{{\k}^{\s}}$). Also, ${\k}^{[\mathbf{1},\ y+z]}$, the centralizer of $[\mathbf{1},\ y+z]$ in $\k$ is equal to ${\sr}^z$. Since ${\mathcal M}\circ{\mathcal V}$ is $K$-equivariant, ${\k}^{{k_0}^{-1}\cdot E'}= {\k}^{[\mathbf{1},\ y+z]}$. Hence, since ${\k}^{E'}= Ad(k_0)({\k}^{{k_0}^{-1}\cdot E'})$, ${\k}^{E'}$ and ${\sr}^z$ are isomorphic Lie algebras. Therefore,
\begin{equation}\label{mainrankeq}
	r_K(\Omega):=\text{rank}\ {\k}^{\Phi(E')} - \text{rank}\ {\k}^{E'}=\text{rank}\ K - \text{rank}\ {\S}^{z}.
\end{equation}
\par
 From Panyushev (Theorem 2.3 in \cite{pa2}) $r_{\Kc}(\mathcal O)=r_{\Kc^x}(\Kc^x\slash \Kc^{\mathfrak s})+ r_{\mathbf S}(Z)$. By Corollary 2(i) of Theorem 1 in \cite{pa1}, $r_{\Kc^x}(\Kc^x\slash \Kc^{\mathfrak s})= \text{rank}\ K^x - \text{rank}\ {\mathbf S}_{\mathbf R}$. By the observation following equation (6.7) in \cite{king}, $r_{\mathbf S}(Z)=\text{rank}\ {\mathbf S}_{\mathbf R} - \text{rank}\ {{\mathbf S}_{\mathbf R}}^{z}$. Hence, by equation(\ref{mainrankeq}),  
\begin{equation}\label{rankequality}
r_{\Kc}(\mathcal O)=\text{rank}\ K - \text{rank}\ {{\mathbf S}_{\mathbf R}}^{z},
\end{equation}
which is the same as $r_K(\Omega)$. \par
	The argument for part (b) of the theorem starts by noting that for all $E'\in\Omega$, Lemma \ref{tanspacedecomp} implies that if $W=T_{E'}(K\cdot E')$, then
\begin{equation}\label{normalsliceeq1}
\dim \Omega=\dim T_{E'}(\Omega)= 2\dim \left({{\k}^{\Phi(E')}\slash {\k}^{E'}}\right) +\dim \left({{\k}\slash {\k}^{\Phi(E')}}\right)+ \dim \frac{W^{\bot}}{W\cap W^{\bot}}.
\end{equation}
Therefore, 
\begin{equation}\label{normalsliceeq2}
\dim \Omega= 2(\dim {\k}^{\Phi(E')}-\dim {\k}^{E'}) +\dim {\k}-\dim {\k}^{\Phi(E')}+ \dim \frac{W^{\bot}}{W\cap W^{\bot}}.
\end{equation}
\par
 Choose $E'\in \Omega_d\cap\Omega_{\Phi}$. By Lemma 2 in section 3 of \cite{huckwurz}, we have \break $[{\k}^{\Phi(E')},\ {\k}^{\Phi(E')}]\subset {\k}^{E'}$. Since ${\k}^{E'}\subset {\k}^{\Phi(E')}$, $[{\k}^{\Phi(E')},\ {\k}^{\Phi(E')}]=[{\k}^{E'},\ {\k}^{E'}]$.\par
	Since ${\k}^{\Phi(E')}$ and ${\k}^{E'}$ have the same maximum semisimple ideal, 
\begin{equation}\label{keydimeq}
\dim {\k}^{\Phi(E')} -\dim {\k}^{E'}=\text{rank}\ {\k}^{\Phi(E')} - \text{rank}\ {\k}^{E'}.
\end{equation}
\par
	Applying Corollary \ref{symsliceeq}, equation (\ref{normalsliceeq2}) becomes
\begin{equation}\label{normalsliceeq3}
\dim \Omega= \dim {\k}^{\Phi(E')}-\dim {\k}^{E'} +\dim {\k}-\dim {\k}^{E'}+ 2\epsilon(\Omega).
\end{equation}
Since $E'\in\Omega_d\cap\Omega_{\Phi}$, equation (\ref{keydimeq}) implies
\begin{equation*}
\dim \Omega= \text{rank}\ {\k}^{\Phi(E')}-\text{rank}\ {\k}^{E'} +\dim {\k}-\dim {\k}^{E'}+ 2\epsilon(\Omega),
\end{equation*}
where $\dim {\k}-\dim {\k}^{E'}$ is the maximal dimension of a $K$ orbit in $\Omega$. 

	 We have just shown that the codimension of the largest $K$ orbit in $\Omega$ is given by the expression
\begin{equation}\label{normalsliceeq6}
\begin{split}
&\text{rank}\ {\k}^{\Phi(E')}-\text{rank}\ {\k}^{E'} + 2\epsilon(\Omega)\\
&= r_{\Kc}(\mathcal O)+ 2\epsilon(\Omega).
\end{split}
\end{equation}
The last assertion follows from part (a). On the other hand, by equation (6.8) in \cite{king}, taking into account equation (\ref{rankequality}) and the fact that $c(\mathcal O)=c_{\Kc^x}(\Kc^x\slash \Kc^{\mathfrak s})+ c_{\mathbf S}(Z)$ (Theorem 2.3 in \cite{pa2}), the codimension of the largest $K$ orbit in $\Omega$ is also given by the expression:
\begin{equation}\label{normalsliceeq7}
r_{\Kc}(\mathcal O)+ 2c(\mathcal O).
\end{equation}
Part (b) of the theorem follows from equating (\ref{normalsliceeq6}) and (\ref{normalsliceeq7}).\par

\end{proof}


\begin{thebibliography}{15}



\bibitem[1]{myky1} I. V. Mykytyuk, \textit{Actions of Borel Subgroups on Homogeneous Spaces of Reductive Complex Lie Groups and Integrability}, preprint, 2001.


\bibitem[2]{king} D. R. King, \textit{Spherical Nilpotent Orbits and the Kostant-Sekiguchi Correspondence}, Trans. Amer. Math Soc., 354(12) (2002), 4909-4920.

\bibitem[3]{karshon} Y. Karshon, \textit{Hamiltonian actions of Lie groups}, thesis, MIT, 1993.

\bibitem[4]{pa1} Dmitri I. Panyushev, \textit{Complexity and rank of homogeneous spaces},  Geometriae Dedicata \textbf{34} (1990), 249--269.

\bibitem[5]{pa2} $\text{\hbox to 1in{\hrulefill}}$, \textit{Complexity and nilpotent orbits} ,  Manuscripta Math. \textbf{83} (1994), 223--237.


\bibitem[6]{huckwurz} A. T. Huckleberry and T. Wurzbacher, \textit{Multiplicity-free complex manifolds}, Math. Ann., 286(1990), 261-280.

\bibitem[7]{gust1} Victor Guillemin and Shlomo Sternberg, \textit{Multiplicity-free spaces}, J. Diff. Geo. 19(1984), 31-56.

\bibitem[8]{gust2} $\text{\hbox to 1in{\hrulefill}}$, \textit{Symplectic techniques in physics}, Cambridge University Press, 1984.


\bibitem[9]{ver} Michele Vergne, \textit{Instantons et correspondance de Kostant-Sekiguchi}, C.R. Acad. Sci. Paris 320(1995), Serie 1, 901--906.


\end{thebibliography}
\end{document}